# The partial fraction representation of 1/(sin(πs/4)2ξ(½ + s)).

# On the Riemann zeta-function, Part III.

By Anthony Csizmazia

E-mail: apcsi2000@yahoo.com

**Abstract**

An odd meromorphic function f(s) is constructed from the Riemann zeta-function evaluated at one-half plus s. The partial fraction expansion, p(s), of f(s) is obtained using the conjunction of the Riemann hypothesis and hypotheses advanced by the author. That compound hypothesis and the expansion p(s) are employed in Part IV to derive the two-sided Laplace transform representation of f(s) on the open vertical strip of all s with real part between zero and four.

**Table of contents.**

**Abstract. Keywords. MSC (Mathematics Subject Classification).**

**Index of symbols.**

**Review of elements of Parts I-II.**

**Introduction**

The conditional implications (*), (**). A general lemma on partial fraction expansions: Absolute convergence, Vanishing at infinity.

**§1 Unconditional results.**

Assumptions and definitions. Bounds on |ζ(z)|⁻¹. A bound for the remainder in a Taylor series. The assumptions for Claims 1.1, 1.2, 2.3 and 2.4. The special cases (′) and (″).

**§2 Conditional results.**

Counterexample. The Hadamard factorization of ξ(½ + s). The monotonicity principle. Conditional theorem 2.1: Vertical vanishing of f(s) – T(s) on $B_i(α)$. Conditional theorem 2.2: Partial fraction representation of f(s).

**Appendix** Complete monotonicity.

**References**

**Index of symbols.**

**Introduction Δ(s).** (Complex domain) S.

**§1**
B(d), z(s), I(s), G′(s), $B_r$(d), w(s), T(s), $B_r$′(d), $F_n(z_0, h, z_0)$, $M(h, z_0, ρ)$, Δ(s, z).

**§2**
$B_i(α)$, k(s), T(s), $B_i$′(α), $Z_i$. $D(t_0, K)$, S′(k).

**Appendix**

E(s). h(z, u, r). ρ(z, u, n). (Half-plane) $H_+$. ρ(z, u).

<h3 style="text-align:center">The partial fraction representation of 1/(sin(πs/4)2ξ(½ + s)).</h3>

<h3 style="text-align:center">On the Riemann zeta-function, Part III.</h3>

**Review of elements of Parts I-II.**

**Part I**

(Complex plane) C.

**§1 Definitions of l(s), a(s), ξ(s), n(s), f(s), b(s), $V(x_0, x_1)$, $V[x_0, x_1]$, V(ε).**

Let s be complex. Define:

$$l(s) := π^{-s/2} sΓ(s/2) = π^{-s/2} \cdot 2Γ(1 + s/2),\ a(s) := l(s)(s-1),\ ξ(s) := (½)a(s)ζ(s),$$

with ζ(s) the Riemann zeta-function,

$$b(s) := \sin(\pi s/4) a(\tfrac{1}{2} + s). \quad n(s) := \sin(\pi s/4) \cdot 2\xi(\tfrac{1}{2} + s). \quad f(s) := 1/n(s).$$

Say $x_0 < x_1$. Let $V(x_0, x_1)$ be the open vertical strip of all s with $x_0 < \mathrm{Re}(s) < x_1$. Define $V[x_0, x_1]$ to be the closed strip of all s with $x_0 \leq \mathrm{Re}(s) \leq x_1$. Set $V(\varepsilon) := V(0, \varepsilon)$ for positive $\varepsilon$.

**The functional equation of $\zeta(s)$. Symmetries of $\zeta(s)$, $\xi(s)$, $n(s)$, $f(s)$.**
$\xi(s)$ is an entire function. $\zeta(s)$ satisfies the functional equation

$$\xi(\tfrac{1}{2} - s) = \xi(\tfrac{1}{2} + s).$$

Thus the entire function $n(s)$ is odd: $n(-s) = n(s)$. Hence the meromorphic function $f(s)$ is odd. $g(s^*) = (g(s))^*$ for g any of $\zeta$, $\xi$, n, f.

## §2 Standard conjectures

### (2.1) The Riemann Hypothesis, RH

In 1859 B. Riemann formulated the following conjecture.
RH: The real part of each nonreal zero of $\zeta(s)$ is one-half.

The Riemann hypothesis has not yet been resolved.

### (2.2) The simple zeros conjecture, SZC.

RH is allied with the unresolved conjecture stated next.

SZC: Each nonreal zero z of $\zeta(s)$ is simple, $\zeta'(z) \neq 0$.

**Definition of $\gamma_n$.** Let $\gamma_1 < \gamma_2 < \ldots < \gamma_n < \gamma_{n+1} \ldots$ enumerate in order of magnitude the distinct imaginary parts $\gamma_n$ of the zeros z of $\zeta(s)$ with $\mathrm{Im}(z) > 0$. If $n \geq 1$, define $\gamma_{-n} = -\gamma_n$.

### (2.3) The Lindelof hypothesis, LH.

The Riemann hypothesis implies the Lindelof hypothesis stated next.

LH: If $\varepsilon > 0$, then $|\zeta(\tfrac{1}{2} + it)| \leq t^\varepsilon$, for large positive t.

LH implies that for $\sigma \geq \tfrac{1}{2}$ and $\varepsilon > 0$: $\zeta(\sigma + it) = O(|t|^\varepsilon)$, for real t with $|t|$ large.

**(2.4)**.

**Definition of $\delta_k$.** $\delta_k := \min\{\gamma_k - \gamma_{k-1}, \gamma_{k+1} - \gamma_k\}$.

## §3

**Definitions of $V_u$, $V_u'$.** Say u is a multiple of four, u = 4w. Let $V_u := V(u, u+4)$. If $u \neq 0, -4$, set $V_u' = V_u$. Let $V_0' := V(\frac{1}{2}, 4)$. Take $V_{-4}' = -V_0'$.

**Definition of the Pochhammer symbol $(z)_n$.** $(z)_n := \Pi_{0 \leq k \leq n-1}(z + k)$, with z complex and the integer $n \geq 1$. Take $(z)_0 = 1$.

**Definitions of $\tilde{c}(4k)$, $c(z)$.** Let k be an integer $\geq 0$.

$$\tilde{c}(4k) := 1/(\pi^{3/4}\Gamma(5/4 + 2k)(2k - \tfrac{1}{4})\zeta(\tfrac{1}{2} + 4k)).$$

$\tilde{c}(4k) > 0$. Define $c(z) := 1/n'(z)$, for z with $n'(z) \neq 0$. $c(4k) = \tilde{c}(4k)(-(\pi^2))^k$. In particular $c(0) = 2^4/(\pi^{3/4}\Gamma(\tfrac{1}{4}) \cdot (-\zeta(\tfrac{1}{2})))$.

**Definition of $P_0(z)$.** Set $P_0(z) := (-1)\sum_{k \geq 1} \tilde{c}(4k)(-(z^2))^k$.

**Definition of the open disk $B(z, r)$.** Say $r \geq 0$. $B(z, r) := \{s: |s - z| \leq r\}$.

**Definition of Z.** Let **Z** be the set of integers.

## §4 Conjectures introduced by the author

The author advances the following conjectures C1-4. They are within the predictions of the GUE model.

**(4.1)**

**Definitions of $\delta_k'$, $I_k(\alpha)$, $S_k(\alpha)$, $T_*(\alpha)$, $t(1)$, $x(t, \alpha)$, $s(t, \alpha)$, $j_k(\alpha)$.**
Assume RH. Set $\delta_k' := \min\{1/\log(|\gamma_k|), \gamma_k - \gamma_{k-1}, \gamma_{k+1} - \gamma_k\}$. Say $0 < \alpha < \frac{1}{2}$. Define $I_k(\alpha)$ to be the open interval $\{t: \gamma_k + \alpha\delta_k' < t < \gamma_{k+1} - \alpha\delta_{k+1}'\}$, $S_k(\alpha)$ to be the semicircle $\{s: |s - i\gamma_k| = \alpha\delta_k', \text{Re}(s) \geq 0\}$ and $T_*(\alpha) := U_{k \geq 1} (I_k(\alpha) \cup S_k(\alpha))$. Specify that $t(1) := \gamma_1 - \alpha\delta_1'$. Say $t > t(1)$. Let $x(t, \alpha)$ be the unique real x with $x + it$ in $T_*(\alpha)$. Set $s(t, \alpha) := x(t, \alpha) + it$. Interpret $\zeta(\tfrac{1}{2} + it)/(t - \gamma_k)$ at $t = \gamma_k$ as $\zeta'(\tfrac{1}{2} + i\gamma_k)$. Let $j_k(\alpha) := \min\{|\zeta(\tfrac{1}{2} + it)/(t - \gamma_k)| : t \text{ is real and } |t - \gamma_k| \leq \alpha\delta_k'\}$. $|\zeta'(\tfrac{1}{2} + i\gamma_k)| \geq j_k(\alpha)$.

**C1 = Conjecture 1** Assume $0 < \alpha < \frac{1}{2}$, and (i), (ii) as follows hold.
(i) There exist $\varepsilon, \lambda$ with $\varepsilon > 0, \lambda > 0$ such that $|\zeta(\tfrac{1}{2} + s(t, \alpha))| > \lambda \cdot t^{-\varepsilon}$, for $t > t(1)$.
**Definition of $\varepsilon_0$.** Let $\varepsilon_0$ be the infimum of such $\varepsilon$.

(ii) $\varepsilon_0 < 3/4$.

**Definition of C′.** Let C′ be the assumption that each of RH, C1 (i) and $\varepsilon_0 < 7/4$ holds.

**(4.2)**

$n'(z)$ and $c(z)$ are even. Hence $n'$ and $c$ are real-valued on the imaginary axis. If RH and SZC hold, then $c(i\gamma_k) = 1/(b(i\gamma_k)\zeta'(½ + i\gamma_k))$.

**Definitions of A and C°.** $A := 2\sum_{k \geq 1} |c(i\gamma_k)|$. $C° := \sum_{k \geq 1} |c(i\gamma_k)|/(\gamma_k^2)$.

**C2 = Conjecture 2**
(i) There exists a real $\varepsilon$ (with $\varepsilon \geq 0$) such that: for any $\sigma > 0$, there is a $K > 0$ with $|\zeta'(½ + i\gamma_k)| > K\gamma_k^{-(\varepsilon + \sigma)}$, for all $k \geq 1$.
**Definition of $\varepsilon_1$.** Let $\varepsilon_1$ be the least $\varepsilon$ as in (i).
(ii) $\varepsilon_1 < 3/4$.

**(4.3)**

**Definition of B°.** $B° := \sum_{k \geq 1} |c(i\gamma_k)|/\delta_k'$.

**Definitions of the partial fraction expansions $p_r(s)$, $p_i(s)$ and $p(s)$.**

$$p_r(s) := (c(0)/s) + \sum_{w \geq 1} c(4w)(1/(s - 4w) + 1/(s + 4w)).$$

Assume C° is finite.

$$p_i(s) := \sum_{k \geq 1} c(i\gamma_k)(1/(s - i\gamma_k) + 1/(s + i\gamma_k)) = 2s\sum_{k \geq 1} c(i\gamma_k)(1/(s^2 + \gamma_k^2)).$$

$$p(s) := \sum_{z: n(z) = 0} (1/n'(z))(1/(s - z)) = p_i(s) + p_r(s).$$

**Definitions of $Z_i$, $Z°$.** Let $Z_i$ be the set of nonreal zeros of $\zeta(½ + s)$. Set $Z° = (4\mathbf{Z}) \cup Z_i$.

**C3 = Conjecture 3**
(i′) There exists a real $\varepsilon$ (with $\varepsilon \geq 0$) such that: for any $\sigma > 0$, there is a $K > 0$ with $|\zeta'(½ + i\gamma_k)| > K\gamma_k^{-(\varepsilon + \sigma)}$, for all $k \geq 1$.
**Definition of $\varepsilon_1$.** Let $\varepsilon_1$ be the least $\varepsilon$ as in (i′).
(i) There exists an $\varepsilon \geq 0$ such that: for any $\sigma > 0$, there is a $K > 0$ with $\delta_k > K\gamma_k^{-(\varepsilon + \sigma)}$, for all $k \geq 1$.
**Definition of $\varepsilon_2$.** Let $\varepsilon_2$ be the least $\varepsilon$ as in (i).
**(ii) $\varepsilon_1 + \varepsilon_2 < 3/4$.**

**(4.4)**

**C4 = Conjecture 4**
There exists an α as in C1 for which (i), (ii) as follows also hold.
(i) There is an $\varepsilon \geq 0$ such that: for any $\sigma > 0$, there is a $K > 0$ such that for all positive integers k: $j_k(\alpha) > K\gamma_k^{-(\varepsilon + \sigma)}$.
**Definition of $\tilde{\varepsilon}_1$.** Let $\tilde{\varepsilon}_1$ be the least ε as in part (i). Assume C3 (i).
(ii) $\tilde{\varepsilon}_1 + \varepsilon_2 \leq 1$.

**Definition of C^.** Let C^ be the compound conjecture that RH, C1, C3 and C4 all hold.

## §5

### Introduction

**Definitions of $\lambda(y)$, $g_0(y)$.** Assume A is finite.
Let y be real. Take $\lambda(y) := 2\sum_{k \geq 1} c(i\gamma_k)\cos(\gamma_k y)$.
Say $y > 0$. Set $g_0(y) := P_0(\pi e^{-2y})$. Let $y < 0$. Define $g_0(y) := \lambda(y) + c(0) - P_0(\pi e^{2y})$.

In Part III, §2, Conditional theorem 2.2 and Part IV, §1, Conditional theorem 1.1 the respective culmination of the proof of each of the following theorems is achieved.

**Conditional theorem 5.1 Partial fraction representation of f(s).** *Assume C^. f(s) = p(s) on C – Z°.*

**Conditional theorem 5.2** *Assume C^. On $V_0$: $f(s) = \int_R d(y)e^{sy}g_0(y)$.*

Now return to Part III.

### Introduction

**Review** Part I, Introduction, definition of f(s), §3; §4: (4.3), (4.4); §5, Introduction.

Let us embark on establishing the conditional partial fraction expansion, p(s), of f(s) of §2, Conditional theorem 2.2 herein, previewed in Conditional theorem 5.1 of Part I, §5, Introduction. The expansion p(s) is discussed in relation to C3 in Part I, §4, (4.3). In Part IV that expansion of f(s) is employed to obtain, in Conditional theorem 1.1, the Laplace transform representation of f(s) on the strip $V_0$. A heuristic derivation of a formal expression for the Laplace density $g_0(y)$ of the latter representation from p(s) is delineated in Part I, §5,

Introduction.

**Review** Part I, §4: (4.2); and (4.4), for C^.

C^ includes C3. C3 implies C2. RH and part (i) of C2 together imply SZC. RH and C2 (i) together imply SZC. Together RH and $|c(i\gamma_k)| < \infty$ for all $k \geq 1$, imply SZC. The Conditional lemma 4.1 of Part I, §4, (4.2), gives that C2 implies $A < \infty$.

**Review** the presentation of $Z°$, $Z_i$ and $p_r(s)$ in Part I, §3. Review Part I, §4, (4.2) and the material on $p(s)$, $p_r(s)$ and $p_i(s)$ in (4.3).

Assume $C° := \sum_{k \geq 1} |c(i\gamma_k)|/(\gamma_k^2)$ is finite. The Conditional claim 4.1 of Part I, §4, (4.3), assure that $p_i(s)$ is analytic except for simple poles at $\pm i\gamma_k$ for $k \geq 1$. Also assume RH. Then SZC holds.

**Definition of Δ(s).** Assume $C°$ is finite and RH. $f(s) - p(s)$ has an analytic extension from $C - Z°$ to $C$. Let $\Delta(s)$ be that extension.

We seek to prove that $f(s) = p(s)$. Each of $\theta = f$ and $\theta = p$ has the symmetries $\theta(s^*) = (\theta(s))^*$ and $\theta(-s) = -\theta(s)$. So one can assume that $Re(s)$ and $Im(s)$ are nonnegative.

In Conditional corollary 2.12 it is established that: C^ implies $\Delta(s) = 0$ for all s.

The entire function $\Delta(s)$ vanishes for all s, if $\Delta(s)$ is bounded and there is a sequence of $s(k)$ such that $\Delta(s(k))$ converges to 0 as $k \to \infty$.

*A fortiori* $\Delta(s) = 0$ for all s, if C^ is assumed and each of the following (*), (**) holds.

**Definition of S.** Let S be the set of z with $|Re(z)| \geq ½$.

**The conditional implications (*), (**).**
(*) Assume $A := \sum_{k \geq 1} |c(i\gamma_k)|$ is finite. Then:
$\Delta(s) \to 0$ uniformly as $|s| \to \infty$ on the domain S.
(**) Assume C^. Then:
$\Delta(s) \to 0$ uniformly as $|Im(s)| \to \infty$ on the critical strip $|Re(s)| < ½$.

A global preview of the path followed to arrive at Conditional corollary 2.12 via (*) and (**) can be attained by reading the relevant definitions and the statements of the results cited next.

We will employ Lemma 1.1 to establish each of the aforementioned cases (*) and (**).

The proof of (*) is achieved in Conditional corollary 2.1 . That proof also relies on Lemma 1.2, Lemma 1.3, Conditional lemma 2.1 and Conditional claim 2.1.

The proof of (**) is obtained in Conditional corollary 2.11. That proof also utilizes Conditional lemma 2.2, Conditional corollary 2.6, Conditional corollary 2.7, Conditional theorem 2.1 and Conditional corollary 2.10.

## §1 Unconditional results.

**Review** Part I, §3.

**Assumptions and definitions of B(d), z(s), I(s) and G′(s).**
Assume Z is an infinite subset of C, with $B(0, r) \cap Z$ finite, for each $r \geq 0$. Let $c_z$ be a nonzero complex number, and $d_z > 0$, for z in Z. Assume that for distinct z, z′ in Z, the disks $B(z, d_z)$ and $B(z', d_{z'})$ are disjoint. Let $B(d) := \bigcup_{z \in Z} B(z, d_z)$. If s is in B(d), take z(s) to be the unique z in Z with s in $B(z, d_z)$. If s is in $C - B(d)$, set I(s) := Z. If s is in B(d), set $I(s) := Z - \{z(s)\}$. Set $G'(s) := \sum_{z \in I(s)} |c_z/(s - z)|$.

**A general lemma on partial fraction expansions.**
**Lemma 1.1** *Say $A' := \sum_{z \in Z} |c_z|/d_z$ is finite.*
(1) **Absolute convergence.** $G'(s) \leq A'$.
(2) **Vanishing at infinity.** *$G'(s)$ converges to 0 uniformly, as $|s|$ grows infinite. Hence so does $\sum_{z \in I(s)} c_z/(s - z)$.*

**Proof of (1).** If z is in I(s), then $|s - z| \geq d_z$.

**Proof of (2).** The dominated convergence theorem (See E.H. Lieb, M. Loss [9]) establishes (2) as follows. Let m be the measure on Z induced by $m(\{z\}) = 1$, for each z in Z. Set $t(s, z) = |c_z/(s - z)|$, if z is in I(s). If s is in B(d), set $t(s, z(s)) = 0$. Then $0 \leq t(s, z) \leq |c_z|/d_z$. Also $\int_Z (|c_z|/d_z) d(m(z)) = A'$. A' is finite. t(s, z) converges pointwise to zero on Z as $|s|$ grows infinite. So (2) holds.

A direct proof of (2) for that specific case is given next. Say $\varepsilon > 0$. Let $r(\varepsilon)$ be the least $r' \geq 0$ such that $\sum_{z \in Z, |z| > r'} |c_z|/d_z < \varepsilon/2$. Then $\sum_{z \in I(s), |z| > r(\varepsilon)} |c_z/(s - z)| < \varepsilon/2$. Assume there is a z in Z with $|z| \leq r(\varepsilon)$. Set $E(s) := \sum_{z \in I(s), |z| \leq r(\varepsilon)} |c_z/(s - z)|$. Then $G'(s) < E(s) + \varepsilon/2$.

Suppose that $|s| \geq R > m(\varepsilon)$, with $m(\varepsilon) = \max\{|z|: z$ is in Z, $|z| \leq r(\varepsilon)\}$. It follows that $E(s) \leq \sum_{z \in Z, |z| \leq r(\varepsilon)} |c_z/(s - z)| \leq (1/(R - m(\varepsilon)))\theta$, with $\theta = \sum_{z \in Z, |z| \leq r(\varepsilon)} |c_z|$. Then $\theta > 0$. Let $R(\varepsilon) := m(\varepsilon) + (2/\varepsilon)\theta$. Take s with $|s| \geq R(\varepsilon)$. Then $E(s) \leq \varepsilon/2$. Thus

G′(s) < ε.

Our present aim is to prove (*) above. See Conditional corollary 2.1.

**Definitions of $B_r(d)$, w(s), T(s) and $B_r′(d)$.**
Fix d with $0 < d \leq 2$. Say w is an integer. Set $B_r(d) := \bigcup_{-\infty < w < \infty} B(4w, d)$. Say s is in $B_r(d)$. Let w(s) be the unique w with s in $B(4w, d)$. Define $T(s) := c(4w(s))/(s - 4w(s))$, for $s \neq 4w(s)$. Each of f(s), $p_r(s)$, and T(s) is odd in s on its domain. Set $B_r′(d) := B_r(d) - 4Z$.

Each of $f(s) - T(s)$ and $p_r(s) - T(s)$ has an analytic extension from $B_r′(d)$ to $B_r(d)$.

**Lemma 1.2** *Each of the following functions converges uniformly to zero as |s| grows infinite.*
(1) *f(s) restricted to S - $B_r(d)$.*
(2) *$p_r(s)$ restricted to C - $B_r(d)$.*
(3) *$p_r(s) - T(s)$ restricted to $B_r(d)$.*

**Proof** We may assume Re(s) is nonnegative.

**Proof of (1)** $f(s) := 1/n(s)$ with $n(s) := b(s)\zeta(½ + s)$. $b(s) := \sin(\pi s/4) \cdot a(½ + s)$ with $a(s) := \pi^{-s/2} \cdot 2\Gamma(1 + s/2)(s-1)$. Restrict s to S - $B_r(d)$.

Claim 1 $1 / |b(s)| \leq |s|^{-7/4}(2e\pi / |s|)^{x/2} K′(d)$.

Proof of Claim 1 Say $z = x + it$ with $x > 0$ and t real. Let $\omega = \arctan(x / |t|)$, for t nonzero, and $\omega = \pi/2$ for $t = 0$. Then $0 < \omega \leq \pi/2$. The Stirling formula for $\Gamma(z)$ gives:

$$|\Gamma(z)| \sim (2\pi/e)^{½} \cdot (|z|/e)^{x - ½} e^{-(\pi/2)|t|} e^{|t|\omega}(1 + O(1/|z|)).$$

(See G. Andrews, R. Askey, R. Roy [2].) Thus for $x \geq x′ > 0$,

$$(e^{(\pi/2)|t|}|\Gamma(z)|)^{-1} \leq (e/|z|)^{x′ - ½}(1 + O(1/|z|)).$$

Subclaim *Assume $d > 0$. Let D be the set of z with $|z - n\pi| \geq d$ for all integers n. There exists an $\varepsilon(d)$ such that for all z in D: $1 / |\sin(z)| \leq e^{-|t|}\varepsilon(d)$, with $t = Im(z)$.*

Proof of Subclaim See Part II, §1. Assume $|t| \geq \delta > 0$.

Then $1 / |\sin(x+it)| = 2 \cdot e^{-|t|}(1 + e^{-2|t|}\varepsilon(x, t))$ with $-1 < \varepsilon(x,t) < 1 + 1/\delta$.

$|\sin(z)|$ has period $\pi$ and is even in z. Also $|\sin(z^*)| = |\sin(z)|$. $|\sin(z)|$ does not

vanish on the compact set K of $z = x + it$ with: $0 \leq x \leq \pi/2$, $0 \leq t \leq \delta$ and $|z| \geq d$. Say K is non empty. $\min\{|\sin(z)|: z \text{ in } K\} > 0$. The Subclaim therefore holds.

Thus Claim 1 is valid.

**Bounds on $|\zeta(z)|^{-1}$.**
(i) There exists an $M > 0$ such that for any $z = x + it$ with $x \geq 1$, t real and $|t| \geq e$, one has $|\zeta(z)|^{-1} < M \cdot (\log|t|)^7$.
(ii) If $x \geq \sigma' > 1$, then $|\zeta(z)|^{-1} \leq \zeta(\sigma')/\zeta(2\sigma')$.

Together Claim 1, (i) and (ii) imply (1) of Lemma 1.2, as detailed next. Say $\sigma > 1/2$. Claim 1 and (ii) imply (1), when $\text{Re}(s) \geq \sigma$. $f(s)$ vanishes with uniform ultra-rapidity, as $\text{Re}(s) \to \infty$. Say s is on $V[½, \sigma]$. Apply Claim 1 and (i). Let t be $\text{Im}(s)$. Assume $0 < p < 7/4$. Then $|f(s)| = O(|t|^{-p})$ for large $|t|$. Thus (1) holds.

Each of (2) and (3) follows from the previous Lemma 1.1.

We now present the principles which will be employed in completing the proof of (*) and then in proving (**) stated above.

**A bound for the remainder in a Taylor series**

Assume that $h(z)$ is analytic on a region $\Omega$. Say $0 < r < \rho$. Suppose that the closed disk of all z with $|z - z_0| \leq \rho$ is a subset of $\Omega$. Say n is a positive integer. $h(z) = (\sum_{0 \leq k \leq n-1} (h^{(k)}(z_0)/(k!)) \cdot (z - z_0)^k) + F_n(z, h, z_0)(z - z_0)^n$, with $F_n(z, h, z_0)$ as follows.

**Definitions of $F_n(z_0, h, z_0)$ and $M(h, z_0, \rho)$.** $F_n(z_0, h, z_0) := h^{(n)}(z_0)/(n!)$. Define $M(h, z_0, \rho) := \max\{|h(z)|: |z - z_0| = \rho\}$.

Say $|z - z_0| \leq r$. Then

$$|F_n(z, h, z_0)| \leq (1/(\rho^{n-1}(\rho - r)))M(h, z_0, \rho).$$

See L. Ahlfors [1].

**The assumptions for Claims 1.1, 1.2, 2.3 and 2.4.** (See below.)

Say $0 < r < \rho < \rho_0$ and $r < \rho' < \rho_0'$. Assume $B(s)$ is analytic on $B(z, \rho_0)$, $B(z) = 0$, $B'(z) \neq 0$ and $B(s) \neq 0$ on $B(z, \rho_0)$ for $s \neq z$. Assume $A(s)$ is analytic and nonzero on $B(z, \rho_0')$. $(AB)'(z) = A(z)B'(z)$. Also $1/A(s)$ is analytic on $B(z, \rho_0')$. Assume that $|s - z| \leq r$.

**Definition of $\Delta(s, z)$.** Set $\Delta(s, z) := (A(s)B(s))^{-1} - ((AB)'(z))^{-1} \cdot (s - z)^{-1}$.

**The special cases (') and (").**

We will later use the following special case (') of the Claims 1.1 – 1.2.

(') $B(s) = \sin(\pi s/4)$, $A(s) = 2\xi(\frac{1}{2} + s)$, $z = 4w$, with $w \geq 1$, $r = d$, with $0 < d \leq 2$, $\rho' = 3$ and $\rho_0 = \rho_0' = 4$. Then $n(s) = A(s)B(s)$ and $f(s) = (A(s)B(s))^{-1}$. Now $|s - 4w| \leq d$ gives $T(s) = ((AB)'(z))^{-1} \cdot (s - z)^{-1}$. Then $f(s) - T(s) = \Delta(s, z)$.

In establishing the Conditional theorem 2.1 we will use the following conditional special case (") in the Conditional claims 2.3 - 2.4.

("): Assume RH and SZC. $B(s) = \zeta(\frac{1}{2} + s)$, $A(s) = b(s)$, $z = i\gamma_k$, with $k \geq 1$, $r = \alpha\delta_k'$, with $\alpha$ as in C1 and C4, $\rho = \rho' = \delta_k'$, $\rho_0 = \delta_k$ and $\rho_0' = \gamma_k$. Then $n(s) = A(s)B(s)$ and $f(s) = (A(s)B(s))^{-1}$. Also $|s - i\gamma_k| \leq \alpha\delta_k'$ gives $T(s) = ((AB)'(z))^{-1} \cdot (s - z)^{-1}$. Then $f(s) - T(s) = \Delta(s, z)$.

**Claim 1.1** $\Delta(s, z) = (A(s))^{-1}((B(s))^{-1} - (B'(z))^{-1} \cdot (s - z)^{-1}) + (B'(z))^{-1} \cdot F_1(s, 1/A, z)$.

**Proof** In $AB - A_0B_0 = A(B - B_0) + (A - A_0)B_0$ replace $A$ with $(A(s))^{-1}$, $B$ with $(B(s))^{-1}$, $A_0$ with $(A(z))^{-1}$ and $B_0$ with $(B'(z))^{-1} \cdot (s - z)^{-1}$. Then

$\Delta(s, z) = (A(s))^{-1}((B(s))^{-1} - (B'(z))^{-1} \cdot (s - z)^{-1}) + ((A(s))^{-1} - (A(z))^{-1}) \cdot (B'(z))^{-1} \cdot (s - z)^{-1}$.

$((A(s))^{-1} - (A(z))^{-1}) \cdot (s - z)^{-1} = F_1(s, 1/A, z)$. So Claim 1.1 is valid.

**Claim 1.2**

$|\Delta(s, z)| \leq (1/|A(s)|) \cdot |(B(s))^{-1} - (B'(z))^{-1}(s - z)^{-1}| + (1/|B'(z)|) M(1/A, z, \rho')/(\rho' - r)$.

**Proof** Apply the bound for the remainder in a Taylor series.

**Eq. (^).**
$$|F_1(s, 1/A, z)| \leq M(1/A, z, \rho')/(\rho' - r).$$

**Lemma 1.3** *Assume $0 < d \leq 2$. Restrict $s$ to $B_r(d)$. $f(s) - T(s)$ converges uniformly to zero as $|Re(s)|$ grows infinite.*

**Proof** We may assume $s = x + it$, with $x, t \geq 0$ and $s$ is not a multiple of 4.

$2\xi(u) = a(u)\zeta(u)$. Let $\sigma = Re(u)$ and $\sigma > 1$. $|\zeta(u)| \geq 2 - \zeta(\sigma)$. $\zeta(\sigma)$ decreases from infinity to one as $\sigma$ increases from 1 to infinity. Let $\sigma_0$ be the unique root thereon of $\zeta(\sigma) = 2$. $1/(2|\xi(u)|)$ is analytic for $Re(u) > \sigma_0$.

Say $4w > \sigma_0 - ½ + 4$. Let $|s - 4w| < d$. Apply Claim (1) to the special case (').

$$f(s) - T(s) = \theta(s)(1/(2\xi(½ + s))) + j^{-1} \cdot F_1(s, 1/(2\xi(½ + s)), 4w),$$

with $j = ((d/ds)(\sin(\pi s/4)))|_{s=4w} = (-1)^w(\pi/4)$, $\theta(s) = 1/\sin(\pi s/4) - j^{-1}(1/(s - 4w))$ and $F_1(s, 1/(2\xi(½ + s)), 4w) := (1/(2\xi(½ + s)) - 1/(2\xi(½ + 4w)))/(s - 4w)$.

Now apply Claim 1.2 to the case (').

$$|f(s) - T(s)| \leq |\theta(s)|(1/(2|\xi(½ + s)|)) + (4/\pi) \cdot |F_1(s, 1/(2\xi(½ + s)), 4w)|.$$

$|\theta(s)|$ has period 4. $\theta(s)$ has an analytic extension from $B(4w, 4) - \{4w\}$ to $B(4w, 4)$, when $w = 0$ and and therefore for all integers $w$. $|\theta(s)| \leq B$ for s on $B_r(d)$, with $B = \max\{|\theta(s)| : |s| \leq 2 \text{ and } Re(s), Im(s) \geq 0\}$.

Say $4w > \sigma_0 - ½ + 4$.

$$|F_1(s, 1/(2\xi(½ + s)), 4w)| \leq \max\{1/(2|\xi(½ + s)|) : |s - 4w| = 3\}.$$

$\sup\{|f(s) - T(s)| : |s - 4w| \leq d\} \leq (B + 4/\pi)m(4w - 3, 3)$. Here, when r is real and $\delta \geq 0$, we define $m(r, \delta) := \sup\{1/(2|\xi(½ + s)|) : Re(s) \geq r, |Im(s)| \leq \delta\}$.

Claim *$m(r, \delta)$ converges ultra-rapidly to zero as $r \to \infty$.*

Proof of Claim Say $r > \sigma_0 - ½$. Then $m(r, \delta) \leq (2 - \zeta(½ + r))^{-1}(â(r, \delta))^{-1}$, with $â(r, \delta) := \inf\{|a(½ + s)| : Re(s) \geq r, |Im(s)| \leq \delta\}$.

Say $s = x + it$, with x large and $|t|$ bounded:

$$|a(½ + s)| = (2^3\pi)^{¼}(x^{7/4})(x/(2\pi e))^{x/2}(1 + \varepsilon/x),$$

with $\varepsilon$ bounded.

**Corollary 1.1** *Restrict s to S. $f(s) - p_r(s)$ converges uniformly to zero, as $|s|$ grows infinite.*

**Proof** Say s is in $B_r(d)$. $f(s) - p_r(s) = f(s) - T(s) - (p_r(s) - T(s))$. $p_r(s) - T(s)$ vanishes uniformly as $|Re(s)|$ grows infinite, by Lemma 1.2 (3). $f(s) - T(s)$ does so by the previous Lemma 1.3.

Say s is in $S - B_r(d)$. Lemma 1.2 (1), (2), yields the asserted convergence of $f(s)$, $p_r(s)$ respectively, and hence of $f(s) - p_r(s)$.

## §2 Conditional results.

**Conditional lemma 2.1** *Assume $A := 2\sum_{k \geq 1} |c(i\gamma_k)|$ is finite. Say $d' > 0$. Restrict s to $|Re(s)| \geq d'$. $p_i(s)$ converges uniformly to zero as $|s|$ grows infinite.*

**Proof** Apply **Lemma 1.1**.

Recall that we assumed $0 < d \leq 2$. (See §1, Definitions of $B_r(d)$, $w(s)$, $T(s)$ and $B_r'(d)$.)

**Conditional claim 2.1** *Assume $A < \infty$. Then $\Delta(s) \to 0$ uniformly as $|s| \to \infty$ on $S - B_r(d)$.*

**Proof** $\Delta(s) = (f(s) - p_r(s)) - p_i(s)$. Implement Lemma 1.2 (1), (2). Also Conditional lemma 2.1 gives $p_i(s) \to 0$ uniformly as $|s| \to \infty$ on $S$.

**Conditional corollary 2.1** *Assume $A$ is finite.*
(*) $\Delta(s) \to 0$ uniformly as $|s| \to \infty$ on $S$.

**Proof** Apply the previous Conditional claim 2.1, if s is on $S - B_r(d)$. Say s is in $B_r(d)$. $\Delta(s) = (f(s) - p_r(s)) - p_i(s)$. The previous Corollary 1.1, Conditional lemma 2.1, assures the convergence of $f(s) - p_r(s)$, $p_i(s)$, respectively, to 0 when $|Re(s)| \to \infty$. Therefore (*) holds.

**Conditional claim 2.2** *Assume RH and A is finite. $\Delta(s)$ is entire. $\Delta(s)$ has the symmetries $\Delta(s^*) = \Delta(s)$ and $\Delta(-s) = -\Delta(s)$. $\lim_{\rho > 0, \rho \to \infty} \Delta(\rho w) = 0$, for any $w \neq \pm i$, with w of unit length.*

**Proof** Together RH and $A < \infty$ imply SZC. So $\Delta(s)$ is entire. Implement the previous Conditional corollary 2.1 to obtain the limit result.

Counterexample Nonetheless, *a priori* it is possible that $\lim_{\rho > 0, \rho \to \infty} |\Delta(\rho v)| = \infty$, for $v = \pm i$. That is elucidated next.

There exist entire functions $E(s)$ of infinite order with $\lim_{\rho > 0, \rho \to \infty} |E(\rho v)| = \infty$ for a unique complex v of unit length, but with $\lim_{\rho > 0, \rho \to \infty} E(\rho w) = 0$ for any w

with |w| = 1 and w ≠ v. See Malmquist [11] and Lindelöf [10].

We now use the E(s) and v of the counterexample to construct an entire function G(s) which shares the properties of Δ established in Conditional claim 2.2, but has $\lim_{\rho > 0, \rho \to \infty}$ (-i)G(ρi) = ∞. Set h(s) := E(sv)·(E(s*v))*. Define G(s) := i(h(-is) – h(is)).

We therefore proceed to prove (**) of the Introduction. (See the Conditional corollary 2.11.)

**Definitions of $B_i(\alpha)$, k(s), T(s), $B_i'(\alpha)$ and $Z_i$.**
Assume 0 < α < ½. Set $B_i(\alpha) := U_{-\infty < k < \infty} B(i\gamma_k, \alpha\delta_k')$. Say s is in $B_i(\alpha)$. Define k(s) to be the unique k such that s is in $B(i\gamma_k, \alpha\delta_k')$. Set $T(s) := c(i\gamma_{k(s)})/(s - i\gamma_{k(s)})$, for s ≠ $i\gamma_{k(s)}$. Let $B_i'(\alpha) := B_i(\alpha) - Z_i$, with $Z_i := \{i\gamma_k: k$ is an integer$\}$. Assume 0 < d ≤ 2. T(s) is defined on $B_r'(d) \cup B_i'(\alpha)$. Note that $B_r(d)$, $B_i(\alpha)$ are disjoint. Hence so are $B_r'(d)$ and $B_i'(\alpha)$. Assume A is finite. $p_i(s) - T(s)$ has an analytic extension from $B_i'(\alpha)$ to $B_i(\alpha)$.

**Review**: See Part I, §4, (4.3), for B°.

**Conditional lemma 2.2** *Assume* $B° := \sum_{k \geq 1} |c(i\gamma_k)|/\delta_k'$ *is finite.*
(1) *Restrict s to V(-½, ½) - $B_i(\alpha)$. $p_i(s)$ converges uniformly to zero as |Im(s)| grows infinite.*
(2) *Restrict s to $B_i(\alpha)$. $p_i(s) - T(s)$ vanishes uniformly as |Im(s)| → ∞.*

**Proof** Apply Lemma 1.1.

**Conditional corollary 2.2** *Assume B° is finite. Restrict s to C - $B_i(\alpha)$.*
(') *$p_i(s)$ converges uniformly to zero as |s| grows infinite.*

**Proof** Say B° < ∞. Then A < ∞ is finite. So (') holds on S. Now apply Conditional lemma 2.2 (1).

Next we develop results to be used to prove (**) stated above.

**The Hadamard factorization of ξ(½ + s).**

The Hadamard factorization of ζ(s) has the form

$$2\xi(s) = \pi^{s/2} \cdot 2^s \cdot e^{(-1 + \gamma/2)s} \cdot \Pi_{z \in Z'} ((1 - s/z)e^{s/z})^{r(z)},$$

with Z' the set of zeros z of ζ(s) in the critical strip V(1) and r(z) the multiplicity of the zero z. See G. Everest, T. Ward [8], p 205, Theorem 9.27. Say p > 1. $\sum_{z \in}$

$_{Z'}$ r(z) / |z|$^p$ is finite. $\xi$(½ + s) = Ke$^{\varphi s}$ · $\Pi_{\theta \varepsilon Z' - ½}$ ((1 - s/θ)e$^{s/\theta}$)$^{r(½ + \theta)}$ for some constants K, φ. So K = $\xi$(½). $\zeta$(x) is nonzero, when 0 < x < 1. $\xi$(½) is nonzero. $\xi$(½ + s) is even in s. Thus $\xi$(½ + s) = $\xi$(½)$\Pi_\theta$ (1 – (s/θ)$^2$)$^{r(½ + \theta)}$. Here θ is in $Z_i$ = Z' - ½ and Im(θ) > 0. Also $\sum_\theta$ r(½ + θ) / |θ|$^p$ is finite.

Assume RH. Then

$$\xi(½ + s) = \xi(½)\Pi_{k \geq 1} (1 + (s/\gamma_k)^2)^{m(k)},$$

with m(k) = r(½ + i$\gamma_k$) and $\sum_{k \geq 1}$ m(k)/($\gamma_k^p$) finite.

**The monotonicity principle.**

Let

$$E(s) = K(\Pi_{1 \leq m \leq N} (1 – s/(i\varphi_m)))\Pi_{k \geq 1} (1 + (s/\theta_k)^2),$$

with: 0 ≤ N ≤ ∞; K ≠ 0; $\varphi_m$ a nonzero real; $\theta_k$ > 0; $\sum_{m \geq 1}$ 1/ |$\varphi_m$| finite, if N = ∞; and $\sum_{k \geq 1}$ 1/$\theta_k^2$ finite. Say s = x + it, with x, t real. Fix t. Set v = x$^2$. Let x ≠ 0.

**Lemma 2.3 Monotonicity principle** |E(s)|, 1 / |E(s)|, *is a strictly increasing, respectively decreasing, function of v.*

**Proof** Say r is real. |1 - s/(ir)|$^2$ = r$^{-2}$(v + (t - r)$^2$).

**Corollary 2.3** *Say a > 0. 1 / |sin((π/a)s)E(s)| strictly decreases as x increases from 0 to a/2.*

**Proof** Apply the monotonicity principle and

$$|\sin(x' + it')|^2 = ½(\cosh(2t') – \cos(2x')),$$

for real x', t'.

**Conditional corollary 2.4** *Assume RH. |$\xi$(½ + s)| is a strictly increasing function of x$^2$, with x = Re(s).*

In the next corollary interpret (s - i$\gamma_k$)f(s) at s = i$\gamma_k$ as 1/n'(i$\gamma_k$).

**Conditional corollary 2.5** *Assume RH. Each of |f(s)| and |(s - i$\gamma_k$)f(s)| is a strictly decreasing function of x$^2$, with x = Re(s) and -2 ≤ x ≤ 2.*

We will use the previous Lemma 2.3, but not the stronger theorem stated next.

The proof of the theorem is in the Appendix.

**Theorem** $|E(x + it)|^{-2}$ is a completely monotone function of $x^2$.

**Review** Part I, §4, (4.1): the asymptotic behavior of $|b(s)|$ on a vertical strip of finite width, Definitions, C1 and C'.

**Conditional lemma 2.4** *Assume C'. Restrict s to $V(-½, ½) - B_i(\alpha)$. $f(s)$ converges uniformly to zero as $|Im(s)|$ grows infinite.*

**Proof** $|f(s)| = |f(-s)| = |f(s^*)|$. So assume $s = x + it$ with $0 \leq x < ½$ and $t > t(1)$. The previous conditional corollary enables us to replace s with $s(t, \alpha)$, the horizontal projection to the left of s onto the contour $T_*(\alpha)$ formed from the union, for $k \geq 1$, of the right-half boundaries $S_k(\alpha)$ of the $B(i\gamma_k, \alpha\delta_k')$ together with the connecting intermediate intervals $I_k(\alpha)$ on the imaginary axis. Implement the previous Conditional corollary 2.5. RH has the consequence that $|f(s)| \leq |f(s(t, \alpha))|$.

$f(z) := (1/b(z))\cdot(1/\zeta(½ + z))$. There exists a positive $K_1$ such that for any x, t satisfying $0 \leq x < ½$ and $t > t(1)$: $|b(x + it)|^{-1} \leq K_1 \cdot |t|^{-7/4}$.

C1 (i) yields that $\varepsilon_0$ is nonnegative and such that for any $\varepsilon > 0$, there is a positive $K(\varepsilon)$ for which $|\zeta(½ + s(t, \alpha))| > K(\varepsilon) \cdot t^p$ for all $t > t(1)$, with $p = \varepsilon_0 + \varepsilon$.

So there is a $K_2(\varepsilon)$ such that for $t > t(1)$: $|f(s(t, \alpha))| < K_2(\varepsilon)\cdot|t|^{-q}$, with $q = (7/4 - \varepsilon_0) - \varepsilon$. Assume $\varepsilon_0 < 7/4$. Take $\varepsilon$ with $0 < \varepsilon < 7/4 - \varepsilon_0$. Then $q > 0$. So the Conditional lemma 2.4 is valid.

**Conditional corollary 2.6** *Assume C'. Restrict s to $C - (B_r(d) U B_i(\alpha))$. $f(s)$ converges uniformly to zero as $|s|$ grows infinite.*

**Proof** Apply Lemma 1.2 (1) and the previous Conditional lemma 2.4.

**Review** Part I, §4, (4.3).

$A < \infty$ implies $C°$ is finite. C3 implies $B°$ is finite. $B° < \infty$ implies A is finite.

**Conditional corollary 2.7** *Assume C' and $B° < \infty$. Restrict s to $C - B_i(\alpha)$. $\Delta(s)$ converges uniformly to zero as $|s|$ grows infinite.*

**Proof** Say s is in S. $B° < \infty$. Then A is finite. So (*) (of the Introduction) holds, by the Conditional corollary 2.1.

$\Delta(s) = f(s) - (p_r(s) + p_i(s))$. Assume s is in $V(-\frac{1}{2},\frac{1}{2})$. Let $|\text{Im}(s)| \to \infty$. $p_r(s)$ vanishes uniformly on $V(-\frac{1}{2},\frac{1}{2})$, by Lemma 1.2 (2). Restrict s to $V(-\frac{1}{2}, \frac{1}{2})$ - $B_i(\alpha)$. $B°$ is finite. That coupled with the Conditional lemma 2.2 (1) entails that $p_i(s)$ vanishes uniformly. C' together with Conditional corollary 2.6 gives $f(s) \to 0$ uniformly. Thus $\Delta(s) \to 0$ uniformly.

In proving that C^ implies $f(s) = p(s)$, the main difficulty is establishing the Conditional theorem 2.1.

On $B_i(\alpha)$ each of $\theta = f$ and $\theta = T$ has the symmetries $\theta(s^*) = (\theta(s))^*$ and $\theta(-s) = -\theta(s)$. So one can assume that Re(s), Im(s) are nonnegative.

The proof of the Conditional theorem 2.1 depends on the series of preliminary results developed next.

The assumptions for the following Conditional claims 2.3 - 2.4 were specified in §1 along with the special case ('') to be considered.

**Conditional claim 2.3** *Assume RH and SZC.*

$$\Delta(s, z) = (B'(z))^{-1} \cdot (-F_2(s, B, z)/(A(s) \cdot (B(s)/(s - z))) + F_1(s, 1/A, z)).$$

**Proof** Apply Claim 1.1.

**Eq.(').**

$$(B(s))^{-1} - (B'(z))^{-1} \cdot (s - z)^{-1} = (1/B'(z)) \cdot (-F_2(s, B, z))/(B(s)/(s - z)).$$

So Claim (1') holds.

**Conditional claim 2.4** *Assume RH and SZC.*

$$|\Delta(s, z)| \leq |B'(z)|^{-1} \cdot (M(B, z, \rho)/(\rho(\rho - r)|A(s) \cdot B(s)/(s - z)|) + M(1/A, z, \rho')/(\rho' - r)).$$

**Proof** Apply Claim 1.2 and Eq.('). Then

$$|\Delta(s, z)| \leq |B'(z)|^{-1} \cdot (|F_2(s, B, z)| / |A(s) \cdot (B(s)/(s - z))| + |F_1(s, 1/A, z)|).$$

Apply the bound for the remainder in a Taylor series to get Eq. (^) and

$$|F_2(s, B, z)| \leq M(B, z, \rho)/(\rho(\rho - r)).$$

Consider the special case (″) specified above in §1.

**Conditional claim 2.5** *Assume RH, SZC and $0 < \alpha < \tfrac{1}{2}$. If $|s - i\gamma_k| \leq \alpha\delta_k'$, then*

$$|f(s) - T(s)| \leq$$

$$(1 - \alpha)^{-1} \cdot |\zeta'(\tfrac{1}{2} + i\gamma_k)|^{-1} \cdot (M(\zeta(\tfrac{1}{2} + s), i\gamma_k, \delta_k')(\delta_k')^{-2}|(s - i\gamma_k)f(s)| + M(1/b(s), i\gamma_k, \delta_k')(\delta_k')^{-1}).$$

**Proof** The current assumptions assure the application of Conditional claim 2.4 to the case (″) specified above.

Next we establish the Corollary 2.9 and Conditional lemma 2.5. They will be used in the proof of the Conditional claim 2.7.

**Claim 2.6** *Say $x_0 \geq 0$ and $\delta > 0$. There exist $K(x_0, \delta)$, $\varepsilon$ such that for any $s = x + it$, with $0 \leq x \leq x_0$, $t$ real and $|t| \geq \delta$:*

$$|\zeta(\tfrac{1}{2} - s^*)| \sim |\zeta(\tfrac{1}{2} + s)| \cdot (|t|/(2\pi))^x (1 + \varepsilon(x, t)/|t|),$$

*with $|\varepsilon(x, t)| \leq K(x_0, \delta)$.*

**Proof** The Claim 2.6 results upon the application of all of the following.
(1) The functional equation $\zeta(1 - u) = 2(2\pi)^{-u}\Gamma(u)\cos(\pi u/2)\zeta(u)$. See Tom M. Apostol [3].
(2) The symmetry $(\zeta(u))^* = \zeta(u^*)$.
(3) The Stirling approximation to $\Gamma(u)$ on a vertical strip of finite width. Say $x_0' \leq x \leq x_0$, $t$ is real and $|t| \geq T > 0$.

$$|\Gamma(x + it)| \sim (2\pi)^{\tfrac{1}{2}} \cdot |t|^{x - \tfrac{1}{2}} \cdot e^{-\tfrac{1}{2}\pi|t|}(1 + \varepsilon(x, t)/|t|),$$

with $|\varepsilon(x, t)| < K(x_0', x_0)$. See G. Andrews, R. Askey, R. Roy [2].
(4) Say $x$, $t$ are real. $|\cos(x + it)| \sim \tfrac{1}{2}e^{|t|}(1 + e^{-2|t|}\varepsilon(x, t))$, with $|\varepsilon(x, t)| \leq 1$.

**Definition of $D(t_0, K)$.**
Say $t_0 > 1$ and $K \geq 0$. Let $D(t_0, K) := \{x + it: t \geq t_0 \text{ and } 0 \leq x \leq K/\log(t)\}$.

**Corollary 2.8** *There exists a $\theta(t_0, K)$ such that $|\zeta(\tfrac{1}{2} - s^*)| \leq \theta(t_0, K)|\zeta(\tfrac{1}{2} + s)|$, for all $s$ in $D(t_0, K)$.*

**Proof** Apply the previous Claim 2.6 with $x_0 = K/\log(t_0)$. $t^x = \exp(x\log(t))$. $x\log(t) \leq K$.

**Definition of S'(k)** Say $k \geq 1$. Let $S'(k)$ be the semi-disk of s with $\mathrm{Re}(s) \geq 0$ and $|s - i\gamma_k| \leq 1/\log(\gamma_k)$.

**Corollary 2.9** *There exists a $\theta$ such that $|\zeta(\frac{1}{2} - s^*)| \leq \theta|\zeta(\frac{1}{2} + s)|$, provided s is in $\cup_{k \geq 1} S'(k)$.*

**Proof** $S'(k) \subseteq [0, 1/\log(\gamma_k)] \times (i \cdot [\gamma_k - 1/\log(\gamma_k), \gamma_k + 1/\log(\gamma_k)]) \subseteq D(e - 1, 2)$.

**Review Part I**, **§2**, **(2.3)**, the Lindelof hypothesis, LH. Part I, §4, (4.1), the Stirling approximation to $\Gamma(z)$.

**Conditional lemma 2.5** *Assume RH. Suppose that: for any $t \geq 1$, $g(t) \geq 0$; and $\lim_{t \to \infty} g(t) = 0$. Given any $\varepsilon > 0$, there exists a $T(\varepsilon) \geq 1$ such that for all x, t with $|t| \geq T(\varepsilon)$ and $0 \leq x \leq g(t)$: $|\zeta(\frac{1}{2} + x + it)| \leq |t|^\varepsilon$.*

**Proof** $|\zeta(u^*)| = |\zeta(u)|$. So assume $t \geq 0$. RH and the Conditional corollary 2.4 together give the following. If $0 \leq x \leq \sigma$, then $|\zeta(\frac{1}{2} + x + it)| \leq |\zeta(\frac{1}{2} + \sigma + it)| \cdot r(\sigma, x, t)$, with $r(\sigma, x, t) := |a(\frac{1}{2} + \sigma + it)/a(\frac{1}{2} + x + it)|$.

Say $\sigma_0 > 0$ and $0 \leq x \leq \sigma \leq \sigma_0$. The Stirling approximation to $\Gamma(z)$ yields that for some $T_0(\sigma_0) > 0$ and all $t \geq T_0(\sigma_0)$: $r(\sigma, x, t) = (t/(2\pi))^{\frac{1}{2}(\sigma - x)} \cdot (1 + \theta(\sigma, x, t)/t)$, with $|\theta(\sigma, x, t)| < K(\sigma_0)$.

Let any $\varepsilon > 0$ be given. Take $T_0(\sigma_0, \varepsilon) \geq 1$ such that for any $t \geq T_0(\sigma_0, \varepsilon)$: $g(t) \leq \sigma$, with $\sigma = \min\{\sigma_0, \varepsilon\}$. If $t \geq T_0(\sigma_0, \varepsilon)$ and $0 \leq x \leq g(t)$, then $|\zeta(\frac{1}{2} + x + it)| \leq |\zeta(\frac{1}{2} + \sigma + it)| \cdot (t^{\varepsilon/2}) \cdot (1 + K(\sigma_0))$.

RH entails LH. Say $0 < \varepsilon' < \varepsilon/2$. LH gives that there is a $T(\sigma, \varepsilon) \geq T_0(\sigma_0, \varepsilon)$ such that for all $t \geq T(\sigma, \varepsilon)$: $|\zeta(\frac{1}{2} + \sigma + it)| \leq t^{\varepsilon'}$. Those t have $|\zeta(\frac{1}{2} + x + it)| \leq t^{\varepsilon/2 + \varepsilon'}(1 + K(\sigma_0))$, when $0 \leq x \leq g(t)$. $\varepsilon/2 + \varepsilon' < \varepsilon$. So the Conditional lemma 2.5 holds.

**Conditional claim 2.7** *Assume RH, SZC and $0 < \alpha < \frac{1}{2}$. If $|s - i\gamma_k| \leq \alpha\delta_k'$ and $\varepsilon > 0$, then there exists a $K'(\varepsilon) > 0$ such that:*

$$|f(s) - T(s)| \leq K'(\varepsilon) \cdot |\zeta'(\frac{1}{2} + i\gamma_k)|^{-1} \cdot (|\gamma_k|^{-(7/4 - \varepsilon)}(\delta_k')^{-2}(j_k(\alpha))^{-1} + |\gamma_k|^{-7/4}(\delta_k')^{-1}).$$

**Proof** The Conditional corollary 2.5 yields that $|(s - i\gamma_k)f(s)| \leq |(t - \gamma_k)f(it)|$, for s in $B(i\gamma_k, \alpha\delta_k')$ and $t = \mathrm{Im}(s)$. It results that $|(s - i\gamma_k)f(s)| \leq (\beta_k)^{-1}(j_k(\alpha))^{-1}$, with $\beta_k := \min\{|b(it)|: t \text{ is real and } |t - \gamma_k| \leq \alpha\delta_k'\}$ and $j_k(\alpha)$ as in Part I, §4, (4.1), Definitions. So the previous Conditional claim 2.5 yields the following. If $|s - i\gamma_k| \leq \alpha\delta_k'$, then

$$|f(s) - T(s)| \leq$$

$(1 - \alpha)^{-1} \cdot |\zeta'(\frac{1}{2} + i\gamma_k)|^{-1} \cdot (M(\zeta(\frac{1}{2} + s), i\gamma_k, \delta_k')(\delta_k')^{-2}(\beta_k)^{-1}(j_k(\alpha))^{-1} + M(1/b(s), i\gamma_k, \delta_k')(\delta_k')^{-1})$.

The asymptotic approximation to $|b(x + it)|$ of Part I, §4, (4.1), has the following consequences. There exists a $K > 0$ such that for any integer k: $\min\{|b(s)|: s \text{ in } B(i\gamma_k, \delta_k')\} \geq K \cdot |\gamma_k|^{7/4}$. Hence $(\beta_k)^{-1}$ and $M(1/b(s), i\gamma_k, \delta_k')$ are each $\sim O(|\gamma_k|^{-7/4})$.

Let θ be as in Corollary 2.9, with θ ≥ 1. Then for all k:

$M(\zeta(\frac{1}{2} + s), i\gamma_k, \delta_k') \leq \theta \cdot \max\{|\zeta(\frac{1}{2} + s)|: |s - i\gamma_k| = \delta_k' \text{ and } Re(s) \geq 0\}$.

The previous Conditional lemma 2.5 has the following special case. Assume RH. Given any ε > 0, there exists a K(ε) such that $|\zeta(\frac{1}{2} + s)| \leq K(\varepsilon)|t|^\varepsilon$, for s in $D(e - 1, 2)$ and $t = Im(s)$.

It follows that RH implies that for any ε > 0, there exists a K(ε) such that for all k: $M(\zeta(\frac{1}{2} + s), i\gamma_k, \delta_k') \leq \theta \cdot K(\varepsilon)|\gamma_k|^\varepsilon$. So the Conditional claim 2.7 is valid.

**Conditional lemma 2.6** *Assume C^. There exist positive σ' and K'' such that for all s in $B_i(\alpha)$: $|f(s) - T(s)| \leq K'' \cdot |s|^{-\sigma'}$.*

**Proof** Assume the compound conjecture C^ that each of RH, C1, C3 and C4 holds (See Part I, §4). C3 subsumes C2 (i). RH and C2 (i) together imply SZC. The previous Conditional claim 2.7 applies.

Say σ > 0. C2 (i) yields that there exists $K_1(\sigma)$ such that for all k: $|\zeta'(\frac{1}{2} + i\gamma_k)|^{-1} \leq K_1(\sigma)|\gamma_k|^r$, with $r = \varepsilon_1 + \sigma$. Also C3 (i) implies there exists $K_2(\sigma)$ such that for all k: $(\delta_k')^{-1} \leq K_2(\sigma)|\gamma_k|^r$, with $r = \varepsilon_2 + \sigma$. C4 (i) gives that there exists $K_3(\sigma)$ such that for all k: $(j_k(\alpha))^{-1} \leq K_3(\sigma)|\gamma_k|^r$, with $r = \tilde{\varepsilon}_1 + \sigma$.

Then for any small positive σ, there exist K(σ), K'(σ) such that for any integer k and any s with $|s - i\gamma_k| \leq \alpha\delta_k'$: $|f(s) - T(s)| \leq K(\sigma) \cdot |\gamma_k|^{-(p-\sigma)} + K'(\sigma) \cdot |\gamma_k|^{-(q-\sigma)}$, with $p = 7/4 - (\varepsilon_1 + 2\varepsilon_2 + \tilde{\varepsilon}_1)$ and $q = 7/4 - (\varepsilon_1 + \varepsilon_2)$. Now C3 (ii) and C4 (ii) together imply that $p > 0$. *A fortiori* C3 (ii) entails that $q > 0$. Take σ with $0 < \sigma < p$. There is a K''(σ) such that for all s in $B_i(\alpha)$: $|f(s) - T(s)| \leq K''(\sigma) \cdot |s|^{-(p-\sigma)}$.

The following theorem is thereby established.

**Conditional theorem 2.1 Vertical vanishing of f(s) − T(s) on $B_i(\alpha)$.** *Assume C^. Restrict s to $B_i(\alpha)$. $f(s) - T(s)$ converges uniformly to zero, as $|Im(s)| \to \infty$.*

**Conditional corollary 2.10** *Assume C^. $\Delta(s) \to 0$ uniformly as $|Im(s)| \to \infty$ on $B_i(\alpha)$.*

**Proof** Restrict s to $B_i(\alpha)$.

$$f(s) - p(s) = (f(s) - T(s)) - (p_i(s) - T(s)) - p_r(s).$$

The previous Conditional theorem 2.1 yields that $f(s) - T(s) \to 0$ uniformly, when $|Im(s)| \to \infty$. $C\hat{}$ includes C3. C3 implies that $B°$ is finite. The Conditional lemma 2.2 (2) gives $p_i(s) - T(s) \to 0$ uniformly as $|Im(s)| \to \infty$. Lemma 1.2 (2) provides that $p_r(s)$ vanishes uniformly as $|Im(s)| \to \infty$.

**Conditional corollary 2.11** *Assume $C\hat{}$.*
*(\*\*) $\Delta(s) \to 0$ uniformly as $|Im(s)| \to \infty$ on the critical strip $|Re(s)| < \frac{1}{2}$.*

**Proof** Say s is on $B_i(\alpha)$. Apply the previous Conditional corollary 2.10.

Restrict s to $V(-\frac{1}{2}, -\frac{1}{2}) - B_i(\alpha)$. $C\hat{}$ implies $C'$ and $B° < \infty$. Implement the Conditional corollary 2.7.

**Conditional corollary 2.12** *Assume $C\hat{}$. $\Delta(s) \to 0$ uniformly as $|s| \to \infty$.*

**Proof** $C\hat{}$ includes C3. C3 entails $A < \infty$. So (\*) holds by the Conditional corollary 2.1. (\*\*) issues from the previous Conditional corollary 2.11.

The next theorem establishes the promised result stated as Conditional theorem 5.1 in Part I, §5, Introduction.

**Conditional theorem 2.2 Partial fraction representation of f(s).** *Assume $C\hat{}$. $f(s) = p(s)$ on $C - Z°$.*

**Proof** $C\hat{}$ implies $C° < \infty$. $C\hat{}$ includes RH. So $\Delta(s)$ is entire. The previous Conditional corollary 2.12 establishes that $\Delta(s)$ is bounded. Therefore $\Delta(s)$ is of constant value. $\Delta(4w + 2) \to 0$, when the integer $w \to \infty$. The latter holds since: $C\hat{}$ implies $A < \infty$; $0 < d \le 2$ and the Conditional claim 2.1 applies. Thus $\Delta(s) = 0$ on C.

**Appendix**

**Complete monotonicity.**

Say: $0 \le N \le \infty$; $K \ne 0$; $\varphi_m$ is a nonzero real; $\theta_k > 0$; $\sum_{m \ge 1} 1/|\varphi_m|$ is finite, if $N = \infty$; and $\sum_{k \ge 1} 1/\theta_k^2$ is finite.

**Definition of E(s).** $E(s) = K(\Pi_{1 \le m \le N} (1 - s/(i\varphi_m)))\Pi_{k \ge 1} (1 + (s/\theta_k)^2)$,

Say s = x + it, with x, t real. Fix t. Set v = x². Let x ≠ 0.

**Theorem** $|E(x + it)|^{-2}$ *is a completely monotone function of* $x^2$.

**Proof** Say n is a positive integer.

**Definitions of h(z, u, r), ρ(z, u, n) and $H_+$.** Say r, u are real, r is nonzero, and Re(z) > 0. Set $h(z, u, r) := r^{-2}(z + (u - r)^2)$. Then $|1 - s/(ir)|^2 = h(v, t, r)$. Take M := min{n, N}.

$$|K(\Pi_{1 \leq m \leq M} (1 - s/(i\varphi_m)))\Pi_{1 \leq k \leq n} (1 + (s/\theta_k)^2)|^{-2} = \rho(v, t, n),$$

with

$$\rho(z, u, n) := |K|^{-2}(\Pi_{1 \leq m \leq M} (h(z, u, \varphi_m))^{-1})\Pi_{1 \leq k \leq n} (h(z, u, \theta_k)h(z, u, -\theta_k))^{-1}.$$

Let $H_+$ be the half-plane of z with Re(z) > 0.

Consider z on $H_+$.

$$1/z = \int_{y > 0} (dy)e^{-zy}.$$

Relative to z, 1/h(z, u, r) is analytic and the Laplace transform relative to y > 0 of a positive function κ(y, u, r). Thus ρ(z, u, n) is the Laplace transform of the positive convolution of the κ(y, u, r) arising from r = $\varphi_m$, ±$\theta_k$, with 1 ≤ m ≤ M and 1 ≤ k ≤ n. Therefore ρ(v, u, n) is completely monotone in v for v positive: $(-1)^j(d/(dv))^j \rho(v, u, n) \geq 0$, for j a nonnegative integer.

The sequence of ρ(z, u, n), with n ≥ 1, converges uniformly on any compact subset of $H_+$.

**Definition of ρ(z, u).**

$$\rho(z, u) := |K|^{-2}(\Pi_{1 \leq m \leq N} (h(z, u, \varphi_m))^{-1})\Pi_{k \geq 1} (h(z, t, \theta_k)h(z, u, -\theta_k))^{-1}.$$

ρ(z, u) is analytic in z on $H_+$.

Let n → ∞. Then $(d/(dz))^j \rho(z, u, n) \to (d/(dz))^j \rho(z, u)$. Hence ρ(v, u) is completely monotone in v for v > 0. $|E(x + it)|^{-2} = \rho(x^2, t)$. So the Theorem is valid.

Say x, t are real. Fix t. Let x ≠ 0.

**Conditional corollary** *Assume RH. $|\xi(½ + x + it)|^{-2}$ is a completely monotone function of $x^2$.*